\documentclass[10pt]{article}
\usepackage{latexsym}\usepackage{amsbsy}\usepackage{amssymb}
\usepackage[latin1]{inputenc}
\usepackage{a4wide}
\newcommand{\eps}{\varepsilon}
\newcommand{\R}{\mathbb R}
\newcommand{\N}{\mathbb N}
\newcommand{\C}{\mathbb C}

\newcommand{\K}{\mathbb K}

\newcommand{\gs}{\ensuremath{{\mathcal G}} }

\newcommand{\es}{\ensuremath{{\mathcal E}} }
\newcommand{\esm}{\ensuremath{{\mathcal E}_M} }

\newcommand{\ns}{\ensuremath{{\mathcal N}} }
\newcommand{\ks}{\ensuremath{{\mathcal K}} }
\newcommand{\rs}{\ensuremath{{\mathcal R}} }
\newcommand{\cs}{\ensuremath{{\mathcal C}} }
\newcommand{\comp}{\subset\subset}
\newcommand{\cinfty}{{\cal C}^\infty}
\newcommand{\al}{\alpha}
 
\newcommand{\ga}{\gamma}
\newcommand{\Ga}{\Gamma}

\newcommand{\vphi}{\varphi}
\newcommand{\intl}{\int\limits}
\newcommand{\D}{{\cal D}}\newcommand{\Vol}{\mbox{Vol\,}}

 \newcommand{\G}{{\cal G}}
 
\newcommand{\CC}{{\cal C}}

\newcommand{\esmh}{\ensuremath{{\mathcal E}^{\mathit h}_M} }

\newcommand{\simh}{\ensuremath{\sim_{\mathit h}}}

\newcommand{\supp}{\mathrm{supp}}

\newtheorem{thr}{\hspace*{-3mm} \bf}[section]
\newcommand{\ethr}{\end{thr}}

\newcommand{\bd}{\begin{thr} {\bf Definition. }}
\newcommand{\ed}{\end{thr}}

\newcommand{\pr}{\noindent{\bf Proof. }}
\newcommand{\proof}{\pr}
\newcommand{\ep}{\hspace*{\fill}$\Box$\ms}

\newcommand{\bthm}{\begin{thr} {\bf Theorem. }}
\newcommand{\ethm}{\end{thr}}

\newcommand{\bp}{\begin{thr} {\bf Proposition. }}
\newcommand{\bc}{\begin{thr} {\bf Corollary. }}
\newcommand{\blem}{\begin{thr} {\bf Lemma. }}
\newcommand{\bex}{\begin{thr} {\bf Example. }\rm} 
\newcommand{\bexs}{\begin{thr} {\bf Examples. }\rm}

\newcommand{\beast}{\begin{eqnarray*}}
\newcommand{\eeast}{\end{eqnarray*}}

\newcommand{\brem}{\begin{thr} {\bf Remark. }\rm}
\newcommand{\ethi}{\end{thr}}

\newcommand{\ms}{\medskip\\}
\newcommand{\nn}{\nonumber}
\newcommand{\beq}{ \begin{equation} }\newcommand{\eeq}{\end{equation} }
\newcommand{\bea}{\begin{eqnarray}}\newcommand{\eea}{\end{eqnarray}}
\newcommand{\beas}{\begin{eqnarray*}}\newcommand{\eeas}{\end{eqnarray*}}
\newcommand{\beqs}{\begin{equation*}}\newcommand{\eeqs}{\end{equation*}}

\newcommand{\ben}{\begin{enumerate}}\newcommand{\een}{\end{enumerate}}
\newcommand{\ba}{\begin{array}}\newcommand{\ea}{\end{array}}
\newcommand{\cd}{{\cal D}}

\begin{document}

\title{Generalized flows and singular ODEs on differentiable manifolds}
\author{Michael Kunzinger \footnote{Electronic mail: michael.kunzinger@univie.ac.at}\\ 
         {\small Department of Mathematics, University of Vienna}\\
         {\small Strudlhofg.\ 4, A-1090 Wien, Austria}\\
        Michael Oberguggenberger \footnote{Electronic mail: michael@mat1.uibk.ac.at}\\
         {\small Department of Engineering Mathematics, Geometry and Computer Science,
           }\\
         {\small University of Innsbruck, Technikerstr.\ 13, A-6020 Innsbruck, Austria}\\ 
        Roland Steinbauer \footnote{Electronic mail: roland.steinbauer@univie.ac.at}\\
         {\small Department of Mathematics, University of Vienna}\\
         {\small Strudlhofg.\ 4, A-1090 Wien, Austria}\\
        James A.\ Vickers \footnote{Electronic mail: J.A.Vickers@maths.soton.ac.uk}\\
         {\small  University of Southampton, Faculty of Mathematical Studies,}\\ 
         {\small Highfield, Southampton SO17 1BJ, United Kingdom}
       }
\maketitle
\date
                                                                           
\begin{abstract}
Based on the concept of manifold valued generalized functions we
initiate a study of nonlinear ordinary differential equations with singular
(in particular: distributional) right hand sides in a global
setting. After  establishing several existence and uniqueness results
for solutions of such equations and flows of singular vector fields we
compare the solution concept employed here with the purely
distributional setting. Finally, we derive criteria securing that
a sequence of smooth flows corresponding to a regularization
of a given singular vector field converges to a measurable limiting flow.
\end{abstract}

\section{Introduction}\label{sec:intro}

The need for considering ordinary differential equations involving
generalized  functions on differentiable manifolds occurs naturally
in a number of applications. Examples include  singular
Hamiltonian mechanics (\cite{marsden}, \cite{marsden2}), symmetry
group analysis of differential equations involving singularities
(\cite{ga},  \cite{DKP}, \cite{symm}) and geodesic equations in
singular space-times in general relativity (\cite{herbertgeo, geo2, 
geo, vickersESI}).

An appropriate setting for developing a theory capable of handling 
this question is provided by Colombeau's theory of nonlinear
generalized functions. Introduced in \cite{c1, c2} primarily as a tool
for the treatment of nonlinear partial differential equations in the
presence of singularities (cf.\ \cite{MObook}), the theory has
undergone a quite substantial ``geometrization'' in recent years owing
to an increasing number of applications in a predominantly geometric
context. For the so-called full version of the construction,
distinguished by the existence of a canonical embedding of the space
of distributions, this restructuring was carried out in \cite{found,  
vim, Jel}. For the special version of the theory, which will form
the basic setting of the present article, a global construction was
developed in \cite{RD, ndg}. 

Finally, an extension of Colombeau's special construction to a
``nonlinear distributional geometry'' capable of
modeling  generalized functions taking values in
differentiable manifolds was given in \cite{gfvm, gprg, gfvm2}. 
A comprehensive presentation of these
developments can be found in \cite{book}. The
aim of the present paper is to extend this setting to a theory of
singular ordinary differential equations on differentiable manifolds.
\medskip

For the convenience of the reader we review the geometric theory of
generalized functions in the following section where we also introduce
our notational conventions. The basic existence theory for singular
ODEs on manifolds and generalized flows is the subject of Section
\ref{sec:ex}. We compare our setting with the purely distributional
framework put forward in \cite{marsden} in Section \ref{sec:D'}.
We introduce a number of notions of association relations for manifold valued
generalized functions in Section \ref{sec:ass} which in turn are 
used in the final Section \ref{sec:lim} to give necessary criteria 
for the limit of a generalized flow to obey the flow property.

\section{Linear and nonlinear distributional geometry}\label{sec:geo}

In this section we collect some basic definitions from linear 
and nonlinear distributional geometry needed in the sequel. 
Our notational conventions will be based on \cite{book}
throughout.
\medskip

In what follows, $C$ will always denote a generic constant.  $X$ will
be a smooth paracompact Hausdorff manifold of dimension $n$. $K
\subset\subset A$ (with $A\subseteq X$) means that $K$ is a compact subset
of $A^\circ$.  For any vector bundle $\pi_X: E \to X$ over $X$ we denote
by $\Gamma(X,E)$ resp.\ $\Gamma_c(X,E)$ the space of smooth (resp.\ 
smooth compactly supported) sections of $E$. The space of differential
operators $\Gamma(X,E) \to \Gamma(X,E)$ is denoted by ${\cal P}(X,E)$
resp.\ ${\cal P}(X)$ in case $E=X\times \R$ (cf.\ \cite{kahn}).  The
space of smooth sections of the tangent bundle $TX$, i.e., the space
of smooth vector fields on $X$ is denoted by $\mathfrak{X}(X)$.  The
volume bundle over $X$ will be written as $\Vol(X)$, its smooth
sections are called one-densities. The space $\D'(X,E)$ of $E$-valued
distributions on $X$ is defined as the dual of the space of compactly
supported sections of the bundle $E^*\otimes\Vol(X)$:
$$
\D'\,(X,E):=[\Ga_c(X,E^*\otimes\Vol(X))]'\,.
$$
\hspace*{0pt}For $E=X\times \R$ we obtain $\D'(X):=\D'(X,E)$, the space of
distributions on $X$.  We have the following isomorphism of
$\CC^\infty(X)$-modules:
\[
\D'(X)\otimes_{\CC^\infty(X)}\Ga(X,E)\, \,\cong\D'(X,E)\,,
\]
i.e., distributional sections may be viewed as sections with
distributional coefficients.
\medskip

Setting $I=(0,1]$ and ${\mathcal E}(X)=\cinfty(X)^I$, we define the
spaces of moderate and negligible nets in ${\mathcal E}(X)$ by \beas
\esm(X)&:=&\{ (u_\eps)_\eps\in{\mathcal E}(X):\ \forall
K\subset\subset X,\ \forall P\in{\cal P}(X)\ \exists N\in\N:
\\&&\hphantom{(u_\eps)_\eps\in{\mathcal E}(X):\ \forall
  K\subset\subset X,\ \forall P\in{\cal P}(\}} \sup_{p\in
  K}|Pu_\eps(p)|=O(\eps^{-N})\}
\\
\ns(X)&:=& \{ (u_\eps)_\eps\in\esm(X):\ \forall K\subset\subset X,\ 
\forall q \in\N_0:\ \sup_{p\in K}|u_\eps(p)|=O(\eps^{q})\}\,.  \eeas

(Note that in the definition of $\ns(X)$ no conditions on the derivatives of 
$(u_\eps)_\eps$ are necessary, cf.\ \cite{ndg}.)
$\gs(X):= \esm(X)/\ns(X)$ is called the (special) Colombeau algebra on
$X$, its elements are written as $u=[(u_\eps)_\eps)]$. $\gs(\_)$ is a
fine sheaf of differential algebras with respect to the Lie derivative
along smooth vector fields.  $\cinfty(X)$ is a subalgebra of $\gs(X)$
and there exist injective sheaf morphisms embedding $\D'(\_)$ into
$\gs(\_)$.

Elements of $\gs(X)$ are uniquely determined by their values on
generalized points in the following way: $(p_\eps)_\eps \in X^I$ is
called compactly supported if $\exists \eps_0$, $K\comp X$ such that
$p_\eps\in K$ for $\eps<\eps_0$; the set of compactly supported points
is denoted by $X_c$. Two nets
$(p_\eps),\, (q_\eps)_\eps\in X_c$ are called equivalent,
$(p_\eps)_\eps\sim (q_\eps)_\eps$, if $d_h(p_\eps,q_\eps) = O(\eps^m)$
for each $m>0$, where $d_h$ denotes the distance function induced on
$X$ by one (hence every) Riemannian metric $h$. The quotient space
$\tilde X_c$ of the set of compactly supported points modulo $\sim$ is
called the space of compactly supported generalized points on $X$ and
we write $\tilde p = [(p_\eps)_\eps]$.  For $X=\R$ we use the notation
$\rs_c$ instead of $\tilde X_c$.  For $\tilde p \in \tilde X_c$, $u\in
\gs(X)$, $[(u_\eps(p_\eps))_\eps]$ gives a well-defined element of
$\ks = \rs$ resp.\ $\cs$ (the space of generalized numbers
(corresponding to $\K=\R$ resp.\ $\C$ and defined as the set of
moderate nets of numbers $(r_\eps)_\eps \in \K^I$ with $|r_\eps| =
O(\eps^{-N})$ for some $N$ modulo negligible nets $|r_\eps| =
O(\eps^{m})$ for each $m$). $u\in \gs(X)$ is uniquely determined by
its point values on $\tilde X_c$, i.e., $u=v\ \Leftrightarrow u(\tilde
p)=v(\tilde p)\ \forall \tilde p\in\tilde X_c$ (\cite{point},
\cite{ndg}).

Colombeau generalized sections of $E\to X$ are defined analogously to $\gs(X)$
using asymptotic estimates with respect to the norm on the fibers of $E$ 
induced by any Riemannian metric $h$ on $X$, which we will denote by $\|\,\|_h$
throughout. Setting $\Gamma_{\mathcal E}(X,E)=\Ga(X,E)^I$ we define 
\beas
        \Gamma_{\esm}(X,E)&:=& \{ (s_\eps)_{\eps}\in \Gamma_{\mathcal E}(X,E) :
                \ \forall P\in {\mathcal P}(X,E)\, \forall K\comp X \, \exists N\in \N:\\
                 &&\hspace{6cm}\sup_{p\in K}\|Pu_\eps(p)\|_h = O(\eps^{-N})\}\\
        \Gamma_\ns(X,E)&:=& \{ (s_\eps)_{\eps}\in \Gamma_{\esm}(X,E) :
                \ \forall K\comp X \, \forall m\in \N:\\
                 &&\hspace{6.5cm}\sup_{p\in K}\|u_\eps(p)\|_h = O(\eps^{m})\}\,.\
\eeas
Then $\Gamma_\gs(X,E):=\Gamma_{\esm}(X,E)/\Gamma_\ns(X,E)$.
$\Ga_\gs(\_\,,E)$ is a fine sheaf of  projective and finitely generated  
$\gs(X)$-modules, moreover
$$ 
 \Gamma_\gs(X,E)=\gs(X)\otimes_{\CC^\infty(X)}\Ga(X,E)\,.
$$
In case $E$ is some tensor bundle $T^r_sX$ we use the notation $\gs^r_s(X)$ for $\Gamma_\gs(X,T^r_sX)$.

Next we turn to the definition of the space of manifold valued generalized functions (\cite{gfvm, gfvm2}). 
Both solutions of generalized ODEs and flows of generalized vector fields will be
modeled as elements of this space. A net $(u_\eps)_\eps \in \CC^\infty(X,Y)^I$ ($Y$ another 
manifold) is called compactly bounded (c-bounded) if for each $K\comp X$ $\exists \eps_0$, $K'\comp Y$
such that $u_\eps(K) \subseteq K'$ for all $\eps<\eps_0$.
The space $\gs[X,Y]$ of c-bounded 
generalized Co\-lom\-beau functions from $X$ to $Y$
is defined as the quotient of the set of $\es_M[X,Y]$ 
of moderate, c-bounded maps from $X$ to $Y$ modulo the equivalence 
relation $\sim$, where
$\esm[X,Y]$ is the set of all $(u_\eps)_\eps \in \cinfty(X,Y)^I$ 
satisfying
\begin{itemize}
  \item[(i)] $(u_\eps)_\eps$ is c-bounded.
  \item[(ii)]
   $\forall k\in\N$, for each chart $(V,\vphi)$ in $X$, each 
   chart $(W,\psi)$ in $Y$, each $L\comp V$ and each $L'\comp W$
   there exists $N\in \N$  with
   $$\sup\limits_{x\in L\cap u_\eps^{-1}(L')} \|D^{(k)}
   (\psi\circ u_\eps \circ \vphi^{-1})(\vphi(p))\| =O(\eps^{-N})\,,
   $$
\end{itemize}
and $(u_\eps)_\eps$ and $(v_\eps)_\eps \in \esm[X,Y]$ are
called equivalent, $(u_\eps)_\eps \sim (v_\eps)_\eps$, if 
\begin{itemize}
\item[(i)] $\forall K\comp X$, $\sup_{p\in K}d_h(u_\eps(p),v_\eps(p)) \to 0$
($\eps\to 0$)
for some (hence every) Riemannian metric $h$ on $Y$.
\item[(ii)] $\forall k\in \N_0\ \forall m\in \N$,
for each chart    $(V,\vphi)$
   in $X$, each chart $(W,\psi)$ in $Y$, each $L\comp V$
   and each $L'\comp W$:
$$
\sup\limits_{x\in L\cap u_\eps^{-1}(L')\cap v_\eps^{-1}(L')}\!\!\!\!\!\!\!\!\!\!\!\!\!\!\!\!\!\!\!
\|D^{(k)}(\psi\circ u_\eps\circ \vphi^{-1}
- \psi\circ v_\eps\circ \vphi^{-1})(\vphi(p))\|
=O(\eps^m).
$$
\end{itemize}
Moderateness and equivalence of nets in $\CC^\infty(X,Y)^I$ can be 
tested equivalently by composition with smooth functions, i.e., 
$(u_\eps)_\eps\in\esm[X,Y]$ iff $(f\circ u_\eps)_\eps\in\esm(X)\ \forall f\in
\CC^\infty(Y)$ (\cite{gfvm2}, 3.2) and two nets $(u_\eps)_\eps$ 
and $(v_\eps)_\eps$ in $\esm[X,Y]$ are equivalent iff $(f\circ u_\eps-f\circ v_\eps)_\eps\in\ns(X)\
\forall f \in\CC^\infty(Y)$ (\cite{gfvm2}, 3.3). 

Inserting a compactly supported point
$\tilde p\in\tilde X_c$ into $u\in\gs[X,Y]$ yields a well-defined element $[u_\eps(p_\eps)_\eps]
\in\tilde Y_c$ and again these generalized point values characterize elements of $\gs[X,Y]$. 
Typically, elements of $\gs[X,Y]$ are capable of modeling jump-discontinuities. 

In order to be
able to form tangent maps of manifold valued generalized functions, the concept of
generalized vector bundle homomorphisms in the following sense is needed (\cite{gfvm}).
Let $\pi_Y: F\to Y$ be a vector bundle over $Y$ and ${\esm}^{\mathrm{VB}}[E,F]$ 
be the set of all $(u_\eps)_\eps$ $\in$ $\mathrm{Hom}(E,F)^I$ satisfying
\begin{itemize}
\item[(i)] $(\underline{u_\eps})_\eps \in \esm[X,Y]$.
\item[(ii)] $\forall k\in \N_0\
\forall (V,\Phi)$
vector bundle chart in $E$,
$\forall (W,\Psi)$ vector bundle chart in $F$,
$\forall L\comp V\
\forall L'\comp W\ \exists N\in \N\ \exists \eps_1>0\
\exists C>0$ with
$$
\|D^{(k)}
(u_{\eps \mathrm{\Psi}\mathrm{\Phi}}^{(2)}(\vphi(p)))\|
\le C\eps^{-N}
$$
for all $\eps<\eps_1$ and all $p\in L\cap\underline{u_\eps}^{-1}(L')$, where $\|\,.\,\|$ 
denotes any matrix norm. 
\end{itemize}
Here, $\underline{u}_\eps$ is the unique element of $\CC^\infty(X,Y)$ such that
$\pi_Y\circ u_\eps = \underline{u}_\eps\circ \pi_X$ and 
$u_{\eps\mathrm{\Psi}\mathrm{\Phi}}:=
\mathrm{\Psi}\circ u_\eps \circ \mathrm{\Phi}^{-1} 
= (x,\xi) \mapsto (u_{\eps\mathrm{\Psi}\mathrm{\Phi}}^{(1)}(x),
u_{\eps\mathrm{\Psi}\mathrm{\Phi}}^{(2)}(x)\cdot\xi)\,.$

$(u_\eps)_\eps$, $(v_\eps)_\eps \in {\esm}^{\mathrm{VB}}[E,F]$
are called $vb$-equivalent, $(u_\eps)_\eps \sim_{vb} (v_\eps)_\eps$, if
\begin{itemize}
\item[(i)] $(\underline{u_\eps})_\eps \sim (\underline{v_\eps})_\eps$ in
$\esm[X,Y]$.
\item[(ii)] $\forall k\in \N_0\ \forall m\in \N\ \forall (V,\Phi)$
vector bundle chart in
$E$, $\forall (W,\Psi)$ vector bundle chart in $F$,
$\forall L\comp V\ \forall L'\comp W
\ \exists \eps_1>0\ \exists C>0$ such that:
$$
\|D^{(k)}(u_{\eps \mathrm{\Psi}\mathrm{\Phi}}^{(2)}
-v_{\eps \mathrm{\Psi}\mathrm{\Phi}}^{(2)})(\vphi(p))\|
\le C\eps^{m}
$$
for all $\eps<\eps_1$ and all $p\in L\cap\underline{u_\eps}^{-1}(L')
\cap\underline{v_\eps}^{-1}(L')$.
\end{itemize}

Then  $\mathrm{Hom}_{\gs}[E,F] := {\esm}^{\mathrm{VB}}[E,F]\big/\sim_{vb}$.
For $u\in \mathrm{Hom}_{\gs}[E,F]$, $\underline{u} :=[(\underline{u}_\eps)_\eps]$
is a well-defined element of $\gs[X,Y]$ uniquely characterized by $\underline{u}
\circ\pi_X = \pi_Y\circ u$. The tangent map $Tu:=[(Tu_\eps)_\eps]$ of any 
$u\in\gs[X,Y]$ is a well-defined element of $\mathrm{Hom}_{\gs}[TX,TY]$. 

Finally, we need the space $\gs^h[X,F]$ of hybrid generalized functions 
defined on $X$ and taking values in the vector bundle $F$ (\cite{gprg});
this space will be used to define the notion of a vector field on a curve.
It is defined as follows: Let $\esmh[X,F]$ the set of
all nets $(u_\eps)_\eps \in \cinfty(X,F)^{(0,1]}$ satisfying (with 
$\underline{u_\eps} := \pi_Y\circ u_\eps$)
\begin{itemize}
\item[(i)] $(\underline{u}_\eps)_\eps$ is c-bounded.
\item[(ii)] $\forall k\in \N_0
\ \forall (V,\vphi)$ chart in $X$ $\forall (W,\Psi)$ vector bundle
chart in $F$ $\forall L\comp V\ \forall L'\comp W
\ \exists N\in \N\ \exists \eps_1>0\ \exists C>0$ such that
$$
\|D^{(k)}(\Psi\circ u_\eps \circ \vphi^{-1})(\vphi(p))\| \le
C\eps^{-N}
$$
for each $\eps<\eps_1$ and each $p\in
L\cap\underline{u_\eps}^{-1}(L')$.
\end{itemize}

In particular, $(u_\eps)_\eps \in \esmh[X,F]$ implies $(\underline{u_\eps})_\eps
\in \esm[X,Y]$.
$(u_\eps)_\eps$, $(v_\eps)_\eps \in \esmh[X,F]$ are called equivalent,
$(u_\eps)_\eps \simh (v_\eps)_\eps$, if the following conditions
are satisfied:
\begin{itemize}
\item[(i)] For each $K\comp X$,
$\sup_{p\in K}d_h(\underline{u_\eps}(p),\underline{v_\eps}(p)) \to 0$ 
for some (hence every) Riemannian metric $h$ on $Y$.
\item[(ii)] $\forall k\in \N_0\ \forall m\in \N\ \forall (V,\vphi)$
chart in $X$, $\forall (W,\Psi)$ vector bundle
chart in $F$, $\forall L\comp V\ \forall L'\comp W
\ \exists \eps_1>0\ \exists C>0$ such that
$$
\|D^{(k)}(\Psi\circ u_\eps\circ\vphi^{-1} - \Psi\circ v_\eps\circ\vphi^{-1})
(\vphi(p))\| \le C\eps^m
$$
for each $\eps<\eps_1$ and each $p\in
L\cap\underline{u_\eps}^{-1}(L')\cap
\underline{v_\eps}^{-1}(L')$.
\end{itemize}
If $u\in \gs[X,Y]$, $v\in \Gamma_\gs[Y,F]$ then $v\circ u := [(v_\eps\circ u_\eps)_\eps]$
is a well-defined element of $\gs^h[X,F]$. We will make use of this fact in analyzing the
flow property of generalized flows (cf.\ Theorem \ref{th:mfflow} below).

\section{Basic existence and uniqueness theorems}\label{sec:ex}

To begin with we consider the system of autonomous nonlinear ODEs on $\R^n$
\begin{equation}
  \label{eq:ode}
  \dot x(t)=F(x(t))
\end{equation}
subject to the initial conditions
\begin{equation}
  \label{eq:ic}
  x(t_0)=x_0\,.
\end{equation}
In contrast to previous treatments in the literature (cf.\ 
\cite{book}, sec.\ 1.5, \cite{HO, Ligeza1, Ligeza2, NedRaj}) we 
seek solutions to (\ref{eq:ode}) in the space $\gs[\R,\R^n]$ of c-bounded generalized functions
(cf.\ \cite{gfvm}, \cite{gprg}, \cite{gfvm2}) rather than in $\gs(\R)^n$. We will therefore 
suppose $F$ to be c-bounded rather than a tempered Colombeau generalized function
to give sense to the composition of generalized functions 
on the right hand side of equation (\ref{eq:ode}).  It is precisely this 
shift in the overall setting which will
allow for the treatment of the flow of a generalized vector field on
a differentiable manifold as a generalized function valued in a smooth manifold. 

To begin with we present basic existence and uniqueness results for the above 
initial value problem. 

\bthm\label{th:locode1}
Let $\tilde x_0\in\rs_c^n$ and let $F = [(F_\eps)_\eps]\in\gs(\R^n)^n$ satisfy
 \begin{itemize}
  \item [(i)] $\exists C, \eps_0 > 0$ such that $|F_\eps(x)|\leq C(1+|x|)$ \ $(x\in \R^n,\, \eps< \eps_0)$, and
  \item [(ii)] $|\nabla F|$ is locally of $L^\infty$-log-type (cf.\ \cite{HO}, Def. 2.3(c)), i.e.,
   \[
     \forall K\comp\R^n\ \sup_{x\in K}|\nabla F_\eps(x)|=O(|\log\eps|)\,.
   \]
 \end{itemize}
 Then the initial value problem (\ref{eq:ode}), (\ref{eq:ic}) has a unique 
 solution in $\gs[\R,\R^n]$. Moreover, $\dot x$ is c-bounded.
\ethm 

\proof We start by establishing {\em existence}. By (i) classical ODE
 theory provides us with globally defined solutions on the level of
 representatives, i.e., for any (fixed) $\eps$ there exists
 $x_\eps\in\CC^\infty(\R,\R^n)$ such that 
 \bea\label{epsode}
  \dot x_\eps(t)&=&F_\eps(x_\eps(t))\\ x_\eps(t_0)&=&x_{0\eps}\nn, 
 \eea 
 where $[(x_{0\eps})_\eps]=\tilde x_0\in\rs^n_c$. From (i) we obtain 
 using Gronwall's lemma $|x_\eps(t)|\leq Ce^{Ct}$; 
 hence $x_\eps$ as well as $\dot x_\eps$ is c-bounded.  
 To show moderateness of $x_\eps$ we write
 \bea\label{eq:mod}
  |\ddot x_\eps(t)|=|\nabla F_\eps(x_\eps(t))|\,|\dot x_\eps(t)|\leq C\eps^{-N}
 \eea
 by the c-boundedness of $x_\eps$ and its derivative and the moderateness of $F$. The
 higher order derivatives of $x_\eps$ are now estimated inductively 
 by differentiating equation (\ref{eq:mod}).

 To prove {\em uniqueness} suppose $[(y_\eps)_\eps]$ is another (c-bounded) 
 solution subject to the same initial conditions. Then since $[(y_{0\eps})_\eps]
 =[(x_{0\eps})_\eps]$ there exist $\tilde n_\eps,n_\eps\in\ns(\R)^n$ such that
 \beas
  (x_\eps-y_\eps)(t)&=
   &x_{0\eps}-y_{0\eps}+\intl_0^t\left(F_\eps(x_\eps(s))-F_\eps(y_\eps(s))
     +n_\eps(s)\right)\,ds\\
   &=&\intl_0^t\tilde n_\eps(s)\,ds
    +\intl_0^t\intl_0^1\nabla F_\eps((1-\sigma)y_\eps(s)+\sigma x_\eps(s))\,d\sigma\,
    (x_\eps-y_\eps)(s)\,ds\,.
 \eeas
 Hence on $|t-t_0|\leq T$ by assumption (ii) for any $m>0$ and $\eps$ sufficiently small we obtain
 \[
  |(x_\eps-y_\eps)(t)|\leq C\eps^m\,e^{-TC\log\eps}\leq C\eps^{m-1}\,.
 \]
\ep

For later use (cf.\ Theorem \ref{th:locflow}) we note the following 
stronger set of conditions that also gives an existence and uniqueness result.

\bc\label{th:locode2}
  Let $\tilde x_0\in\rs_c^n$, $F = [(F_\eps)_\eps]\in\gs(\R^n)^n$ and suppose that there exist
 $C,\, \eps_0 > 0$ such that
 \begin{itemize}
  \item [(i)]  $|F_\eps(0)|\leq C$ \ $(\eps< \eps_0)$, and
  \item [(ii)] $|\nabla F_\eps(x)|\leq C$ \ $(x\in \R^n,\, \eps< \eps_0)$.
 \end{itemize}
 Then the initial value problem (\ref{eq:ode}), (\ref{eq:ic}) has a unique 
 solution in $\gs[\R,\R^n]$. Moreover, $\dot x$ is c-bounded.
\ethm

\proof From (i) and (ii) it follows that $F$ in fact satisfies the hypotheses of 
Theorem \ref{th:locode1}.
\ep

Next we give the basic theorem on the flow of system (\ref{eq:ode}).

\bthm\label{th:locflow}
 Let $F\in\gs(\R^n)^n$ satisfy the assumptions (i) and (ii) 
 of Theorem \ref{th:locode1}. 
 Then there exists a unique generalized function $\Phi\in\G[\R^{n+1},\R^n]$,
 the {\em generalized flow} of system (\ref{eq:ode}) such that
 \bea
  \frac{d}{dt}\Phi(t,x)&=&F(\Phi(t,x))\quad\mbox{ in }\G[\R^{1+n},\R^n]\label{eq:flow1}\\
  \Phi(0,.)&=&\mathrm{id}_{\R^n} \quad\mbox{ in }\G[\R^n,\R^n]   \label{eq:flow2}\\
  \Phi(t+s,.)&=&\Phi(t,\Phi(s,.)) \quad\mbox{ in }\G[\R^{2+n},\R^n]\,. \label{eq:flow3}
 \eea
 Moreover, $\frac{d}{dt}\Phi$ is c-bounded and under the assumptions of 
 Corollary \ref{th:locode2} $\nabla_x\Phi$ is c-bounded as well.

\ethm

 As usual we shall often write $\Phi_t$ instead of $\Phi(t,.)$ and use the
 notation $\Phi_t=[(\Phi^\eps_t)_\eps]$.
\vskip6pt
\proof
 Classical theory provides us with a unique and globally defined flow
 $\Phi^\eps$ for fixed $\eps$.
 
 To prove {\em existence}, we conclude from the integral equation corresponding 
 to (\ref{eq:flow1}) that  $\Phi^\eps$ and $\frac{d}{dt}
 \Phi^\eps$ are c-bounded as functions in $(t,x)$. The higher 
 order $t$-derivatives are estimated as in the proof of Theorem 
 \ref{th:locode1}. To estimate the $x$-derivatives we write
 \beq
  \nabla_x\Phi^\eps(t,x)= x+\intl_0^t\nabla F_\eps(\Phi^\eps(s,x))
                         \nabla_x\Phi^\eps(s,x)\,ds\,. \label{eq:xder}
 \eeq 
 Since $\Phi^\eps$ is c-bounded, (ii) and Gronwall's inequality imply on any $
 \tilde K=[0,T]\times
 K\comp\R^{1+n}$
 \[
  |\nabla_x\Phi(t,x)|\leq Ce^{-CT\log\eps}=O(1/\eps^{CT})\,.
 \]
 For $F$ satisfying the assumptions of Corollary \ref{th:locode2}, an analogous estimate
 establishes c-boundedness of $\nabla_x \Phi$.

 The higher order $x$-derivatives are now estimated by successively 
 differentiating equation (\ref{eq:xder}) and using the estimates already obtained. 
 Similarly, the mixed $x,t$-derivatives may be estimated by differentiating 
 the equations for the $x$-derivatives with respect to $t$.
  
  
 To prove {\em uniqueness} assume that $\Psi$ is another solution in
 $\gs[\R^{1+n},\R^n]$. Then fixing any $\tilde x_0=[(x_0)_\eps] \in \rs^n_c$, 
 both $t\mapsto\Phi(t,\tilde x_0)$ and $t\mapsto\Psi(t,\tilde x_0)$) solve 
 the initial value problem
 \beas
  \dot x(t)&=&F(x(t))\\
  x(0)&=&\tilde x_0\,.
 \eeas
 By the uniqueness part of Theorem \ref{th:locode1} 
 we have for all $\tilde x\in\rs_c^n:\ \Phi(.,\tilde x)=\Psi(.,\tilde x)$
 in $\gs[\R,\R^n]$. Hence by \cite{gfvm2}, Th.\ 3.5, 
 $\Phi(\tilde t, \tilde x)=\Psi(\tilde t,\tilde x)$
 for all $(\tilde t, \tilde x) \in\rs_c^{1+n}$.
 Therefore, another appeal to  \cite{gfvm2}, Th.\ 3.5 establishes 
 $\Phi=\Psi$ in $\gs[\R^{1+n},\R^n]$.
  
 Finally, the {\em flow properties} (\ref{eq:flow2}), (\ref{eq:flow3}) hold 
 on the level of representatives by the classical theory. Hence again by the 
 point value characterization \cite{gfvm2}, Th.\ 3.5, the claim follows.
\ep

In order to prove analogous theorems on a manifold we 
introduce the following notions of boundedness in terms of Riemannian metrics on $X$.

\bd Let $\xi\in\gs^1_0(X)$.
\begin{itemize}
\item [(i)] We say that $\xi$ is {\em locally bounded} 
  resp.\ locally of {\em $L^\infty$-log-type} 
  if for all $K\comp X$ and one (hence every) Riemannian metric $h$ on $X$
  we have for one (hence every) representative $\xi_\eps$ 
  \[
   \sup_{p\in K}\|\,\xi_\eps|_p\,\|_h\,\leq\,C\, \quad \mbox{resp.} 
   \quad \sup_{p\in K}\|\,\xi_\eps|_p\,\|_h\,\le C |\log \eps|\,,
  \]
  where $\|\quad\|_h$ denotes the norm induced on $T_pX$ by $h$. 
\item [(ii)] $\xi$ is called 
  {\em globally bounded} with respect to $h$ if for some (hence every) representative 
  $(\xi_\eps)_\eps$ of $\xi$ there exists $C>0$ with
  \[
   \sup_{p\in X}\|\,\xi_\eps|_p\,\|_h\,\leq\, C \,.
  \]
\end{itemize}
\ed
Contrary to the local notions in (i) above, global boundedness obviously depends on 
the Riemannian metric $h$.

\bthm\label{th:mfode}
  Let $(X,h)$ be a complete Riemannian manifold, $\tilde x_0\in \tilde X_c$ 
  and $\xi\in\gs^1_0(X)$ such that
  \begin{itemize}
  \item[(i)] $\xi$ is globally bounded with respect to $h$.
  \item[(ii)] For each differential operator $P\in{\cal P}(X,TX)$ of first order $P\xi$ 
  is locally of $L^\infty$-log-type.
  \end{itemize}
  Then the initial value problem
\bea\label{eq:mfode}
   \dot x(t)&=&\xi(x(t)) \\
   x(t_0)&=&\tilde x_0\nn
\eea

  has a unique solution $x$ in $\gs[\R,X]$.
\ethm

Note that equality (\ref{eq:mfode}) holds in the space ${\mathfrak X}_{\gs}(x)$
of generalized sections along the generalized mapping $x\in\gs[\R,X]$,
defined by (cf. \cite{gprg} Def.\ 4.6) ${\mathfrak X}_{\gs}(x):=
\{v\in\gs^h[\R,TX]\,|\,\pi_X\circ v=x \}$.
\medskip

\pr Choose a representative $(\xi_\eps)_\eps$ of $\xi$. By (i) each $\xi_\eps$
is globally bounded with respect to $h$. Then due to the completeness of $(X,h)$, 
for each $\eps\in I$ there exists a globally defined solution $x_\eps$ of
\begin{equation} \label{odeeps}
  \begin{array}{rcl}
   \dot x_\eps(t)&=&\xi_\eps(x_\eps(t)) \\[3pt]
   x_\eps(t_0)&=&\tilde x_{0\eps}
  \end{array}
\end{equation}
(cf.\ \cite{michor}, Ch.\ 5, R20).

Let $t_1 < t_2 \in \R$. Then denoting by $L$ the length of a curve we have from (i)
\begin{equation} \label{length}
L(x_\eps|_{[t_1,t_2]}) = \int_{t_1}^{t_2} \|\dot x_\eps(s)\|_h\, ds 
= \int_{t_1}^{t_2} \|\xi_\eps(x_\eps(s))\|_h\, ds \le C |t_2 - t_1|
\end{equation}
for all $\eps$. Let $K\comp X$, $\eps_0>0$ such that $x_\eps(t_1) \in K$ for all
$\eps < \eps_0$. Then by the above $\bigcup_{\eps<\eps_0} x_\eps[t_1,t_2]
\subseteq \{p\in X \mid d_h(p,K) \le C|t_2-t_1|\}$. 
Since the latter set is compact by the Hopf-Rinow theorem, it follows that
$(x_\eps)_\eps$ is c-bounded. Due to this fact, moderateness of $(x_\eps)_\eps$ follows 
as in the local case taking into account \cite{gfvm}, Def.\ 2.2, 
using the moderateness of $(\xi_\eps)_\eps$ and applying the differential 
equation for $x_\eps$ inductively.

To establish {\em uniqueness}, let $a>0$ and choose 
$\eps_0 \in I$, $K \comp X$ such that 
$x_\eps([-a-1,a+1]) \cup y_\eps([-a-1,a+1]) \subseteq K$ for all $\eps < \eps_0$.
Let $t_0\in (-a,a)$ and suppose that $(x_\eps)_\eps$ satisfies (\ref{odeeps}), 
and $(y_\eps)_\eps \in \esm[\R,X]$ solves
\begin{equation} \label{odeeps2}
  \begin{array}{rcl}
   \dot y_\eps(t)&=&\xi_\eps(y_\eps(t)) + n_\eps(t) \\[3pt]
   y_\eps(t_0)&=&\tilde y_{0\eps}
  \end{array}
\end{equation}
Here, $\pi_X\circ n_\eps = y_\eps$ for each $\eps$ and 
$(n_\eps)_\eps\sim_h (0\circ y_\eps)_\eps$ in $\esm^h[\R,TX]$, 
where $0$ denotes the zero element in $\Gamma_\gs(X,TX)$ (cf.\ 
\cite{gfvm2}, Prop.\ 5.7).
Also $(\tilde x_{0\eps})_\eps$, $(\tilde y_{0\eps})_\eps 
\in X_c$ satisfy $[(\tilde x_{0\eps})] = [(\tilde y_{0\eps})]$ in $\tilde X_c$.
By \cite{aubin}, Th.\ 1.36 there exists some $r>0$ such that
$K$ can be covered by finitely many metric balls $B_r(p_i)$ $(p_i\in K, 1\le i \le k)$ with
each $B_{4r}(p_i)$ a geodesically convex domain for the chart $\psi_i:=\exp_{p_i}^{-1}$. 
Choose $\eps_1<\eps_0$ such that 
$$
d_h(x_\eps(t_0),y_\eps(t_0)) = d_h(x_{0\eps},y_{0\eps}) < r 
$$
for all $\eps<\eps_1$. 
With $C$ as in (\ref{length}) we choose $0<d< 
\min(r/C,1)$. For fixed $\eps < \eps_1$ there exists some  $i\in \{1,\dots,k\}$ 
(depending on $\eps$) with $x_{0\eps} \in
B_r(p_i)$. Then by convexity, for each $t$ with $|t-t_0| < d$ the entire line connecting 
$\psi_i(x_\eps(t))$ and $\psi_i(y_\eps(t))$ is contained in $\psi_i(B_{3r}(p_i))$.
Given any $m>0$, we may therefore employ the Gronwall argument from the proof of 
   Theorem \ref{th:locode1} to conclude that there exists 
$\eps_2<\eps_1$ such that for $\eps < \eps_2$ 
$$
|\psi_i\circ x_\eps(t) - \psi_i\circ y_\eps(t)| \le C' \eps^m\,.
$$ 
Here, $\eps_2$, $C'$ only depend on $n=[(n_\eps)_\eps]$, $K$, 
$(\tilde x_{0\eps})_\eps$, $(\tilde y_{0\eps})_\eps$ and $\psi_i$ (on a 
compact subset of its domain), hence can be chosen uniformly in $i\in \{1,\dots,k\}$ 
and $t\in [t_0-d,t_0+d]$. Therefore, 
$$
\sup_{t\in [t_0-d,t_0+d]} d_h(x_\eps(t),y_\eps(t)) \le C'' \eps^m\,,
$$
for $\eps<\eps_2$, so $(x_\eps)_\eps \sim (y_\eps)_\eps$  on
$(t_0-d,t_0+d)$ by \cite{gfvm2}, Th.\ 3.3. Since $d$ depends exclusively on $K$ and $C$ it follows
that if $x$ and $y$ coincide in any $t_0\in (-a,a)$ then in fact they agree on an interval 
of fixed minimal length around $t_0$, hence they are identical on all of $(-a,a)$. Since $a$ was
arbitrary it follows that $x$ and $y$ agree globally as elements of $\gs[\R,X]$.
\ep
Based on this result we are now able to establish the following flow
theorem in the global context.
\bthm\label{th:mfflow}
  Let $(X,h)$ be a complete Riemannian manifold and suppose that $\xi\in\gs^1_0(X)$ satisfies
  conditions (i) and (ii) of Theorem \ref{th:mfode}.
  Then there exists a unique generalized function $\Phi\in\gs[\R\times X,X]$,
  the {\em generalized flow of $\xi$}, such that

  \bea\label{eq:mfflow1}
   \frac{d}{dt}\Phi(t,x)&=&\xi(\Phi(t,x))\quad\mbox{ in }\G^h[\R\times X,TX] \\
   \Phi(0,.)&=&\mathrm{id}_{X} \quad\mbox{ in }\G[X,X]   \label{eq:mfflow2}\\
   \Phi(t+s,.)&=&\Phi(t,\Phi(s,.)) \quad\mbox{ in }\G[\R^2\times X,X]\,. 
                                                         \label{eq:mfflow3}
  \eea 
\ethm
\pr {\em Existence:} Choosing a representative $(\xi_\eps)_\eps$ such that each
$\xi_\eps$ is globally bounded with respect to $h$ we obtain a global smooth
flow $\Phi^ \eps$ for each $\eps$. (\ref{eq:mfflow1})--(\ref{eq:mfflow3}) then clearly
hold componentwise for $(\Phi^\eps)_\eps$. Since $h$ is complete, an argument as in
(\ref{length}) shows that any compact subset of $\R\times X$ remains bounded
(hence relatively compact) upon application of $\Phi^\eps$, uniformly
in $\eps$. Thus $(\Phi^\eps)_\eps$ is c-bounded.

To show moderateness of $(\Phi^\eps)_\eps$, we first note that $t$-derivatives
of $\Phi^\eps$ may be estimated according to \cite{gfvm}, Def.\ 2.2,
precisely as in the proof of Theorem \ref{th:mfode}. Next, let
$[0,t']\comp \R$, $K\comp X$ be given and fix $p\in K$. Then there
exist $t_0=0$, $t_1$, \dots, $t_k=t'$ such that each
$\{\Phi^\eps(t,p)\mid t_i \le t \le t_{i+1}\}$ lies entirely within a
chart domain. We may therefore iterate an integral argument as
in (\ref{eq:xder}) to obtain a moderateness estimate on the
first (local) $x$-derivative of $\Phi^\eps$ on $[0,t']\times
\{p\}$. Since only finitely many charts are needed to cover
$[0,t']\times K$ and the constants in the resulting estimates can be
chosen uniformly in $p\in K$, we obtain the moderateness estimate for
first order $x$-derivatives of $\Phi^\eps$ (the case $t'<0$ is treated
analogously). Higher order $x$-derivatives as well as mixed $x$, $t$
derivatives are estimated in the same manner, so $(\Phi^\eps)_\eps$
is indeed moderate. Moreover, we conclude from the above that $(\frac{d}{dt} 
\Phi^\eps)_\eps \in \es_M^h[\R\times X, TX]$. Also, the composition on the right hand
side of (\ref{eq:mfflow1}) yields a well-defined element of
$\gs^h[\R\times X,TX]$ by \cite{gprg}, Th.\ 4.2. (\ref{eq:mfflow1})
therefore holds since it was already established on the level of representatives. 
Similarly, (\ref{eq:mfflow2}), (\ref{eq:mfflow3}) hold for $\Phi =
[(\Phi^\eps)_\eps]$.

Finally, {\em uniqueness} of the flow follows from the point value
characterization of manifold valued generalized functions and Theorem 
\ref{th:mfode}, precisely as in the proof of \ref{th:locflow}. 
\ep

\bd\label{def:g-complete}
 We call a generalized vector field $\xi\in\gs^1_0(X)$ {\em $\gs$-complete}
 if there exists a unique global generalized flow $\Phi\in\gs[\R\times X,X]$ 
 satisfying (\ref{eq:mfflow1}), (\ref{eq:mfflow2}), (\ref{eq:mfflow3}).
\ed
\section{The distributional setting} \label{distsect}\label{sec:D'}

Our next aim is an analysis of the interrelation between
the theory introduced in the previous section and a
purely distributional approach, as provided by Marsden in \cite{marsden}. 

Any distributional theory of ordinary differential equations 
on manifolds faces a number of principal obstacles resulting 
from the basic structure of the theory of
distributions itself. In fact, consider the initial value problem
\begin{equation} \label{eq:distode}
  \begin{array}{rcl}
   \dot x(t)&=&\zeta(x(t)) \\[3pt]
   x(t_0)&=& x_0
  \end{array}
\end{equation}
with $\zeta \in \cd'(X,TX)$ a distributional vector field. The first question
to be answered in treating this problem is in which setting the solution $x$ is
to be sought (there is no concept of distributions taking values in a differentiable
manifold). A similar problem occurs upon trying to introduce a notion of distributional
flow for (\ref{eq:distode}).
Marsden in \cite{marsden} employs a regularization approach to cope
with these problems,  introducing a sequence of smooth vector fields
$\xi_\eps$ approximating $\zeta$. Each $\xi_\eps$ has a classical flow
$\Phi^\eps$ and under certain assumptions the assignment 
$\Psi=\lim_{\eps\to 0}\Phi^\eps$ allows one to associate
a measurable function $\Psi$ to the distributional vector field
$\zeta$.  
However, the question arises under which conditions on $\zeta$, resp.\ the
regularizing sequence, the limiting map $\Psi$ is indeed  a flow, i.e., satisfies
$\Psi_{t+s}=\Psi_t\circ \Psi_s$. The answer provided by Th.\ 6.2 in
\cite{marsden} turns out to be wrong as we shall see
below by an explicit counter-example. This fact is particularly unfortunate as
the main flow theorems in \cite{marsden} both in the general (Th.\ 6.3) and
in the Hamiltonian case (Th.\ 8.4) rest upon Th.\ 6.2. 

With a view to a smooth presentation of these considerations we first recall
the following definition (\cite{marsden}, Def.\ 6.1, with the index set of the 
regularizing sequence changed from $\N$ to $I$ to ease
comparison with the present setting):


Let $\zeta\in\D'(X,TX)$ be a distributional vector field on the manifold $X$
and let $(\xi_\eps)_\eps$ be a net of smooth vector fields with complete flows 
$\Phi^\eps(t,.)$ and $\xi_\eps\to\zeta\in\D'(X,TX)$. 
$\zeta$  is called a {\em vector field with measurable flow} $\Psi_t$ if
\begin{itemize}
 \item [(i)] $\Phi^\eps(t,.)\to\Psi(t,.)$ almost everywhere on $X$ 
             for all $t$ (in particular, $\Psi_t$ is measurable), and
 \item [(ii)] For each $t\in \R$ and each $C\comp X$ there exists
   $\eps_0\in I$ and $K\comp X$ with $C\subseteq K$ such that 
   $\Phi^\eps(t,C)\subseteq K$ for all $\eps$.
\end{itemize}
Note that in our terminology, (ii) says that $\Phi^\eps(t,\,.\,)$ is
c-bounded. Moreover, if $(\xi_\eps)_\eps$ is additionally supposed to
be moderate and $\xi = [(\xi_\eps)_\eps]\in \gs^1_0(X)$ then
$\xi_\eps\to\zeta\in\D'(X,TX)$ is the same as requiring that $\xi$ is
associated with $\zeta$ (see Section \ref{sec:ass} below). 
As remarked in \cite{marsden}, $\Psi$ in
general depends on the chosen regularizing net $(\xi_\eps)_\eps$.

The basic theorem on flows of distributional vector fields then takes
the following form

{\bf Theorem 6.2 of \cite{marsden}}\\
{\em Let $\zeta\in\D'(X,TX)$ be a vector field with measurable flow $\Psi_t$; then
the flow property holds in the following sense
\[
 \Psi_{t+s}=\Psi_t\circ \Psi_s \mbox{ almost everywhere on $X$, }\forall s,t\in\R.
\]  
}
\bigskip

In order to analyze the validity of this claim we consider the
following initial value problem on $X=S^1$
\bea\label{Heq}
 \dot x(t)&=&\zeta(x(t))\\
 x(0)&=&e^{i\alpha_0},\label{Hiv}
\eea
with the vector field $\zeta$ given by
\beq\label{zeta}
 \zeta(e^{i\alpha})=\left(e^{i\alpha},H(\alpha+\frac{\pi}{2})-
                               H(\alpha-\frac{\pi}{2})\right),
\eeq
where $H$ denotes the Heaviside function.

We proceed by replacing $H$ by a suitable regularization.
We choose a {\em scaling function} $\sigma:(0,\infty)\to (0,\infty)$ 
satisfying $\sigma(\eps)\to 0\ (\eps\to 0)$ and a mollifier $\rho\in\D(\R)$ with 
$\rho\geq 0$, $\supp(\rho)\subseteq [-1,1]$ and $\int\rho=1$.
Then we set
\bea
 \rho_{\sigma(\eps)}&:=&\frac{1}{\sigma(\eps)}\,\rho\left(\frac{x}{\sigma(\eps)}\right)\ 
 \mbox { and finally}\\
 H_\eps(x)&:=&\int\limits_{-\infty}^x\rho_{\sigma(\eps)}(s)\,ds.
\eea
Equipping $S^1$ with the standard metric we have the following

\bp\label{counterex} 
Let $\alpha_0\in[-\pi,\pi]$. The initial value problem 
\bea\label{Hepseq}
 \dot x(t)&=&\xi(x(t))\\
 x(0)&=&e^{i\alpha_0},\label{Hepsiv}
\eea
with the vector field $\xi=[(\xi_\eps)_\eps]$ given by
\beq\label{xieps}
 \xi_\eps(e^{i\alpha})=\left(e^{i\alpha},H_\eps(\alpha+\frac{\pi}{2})-
                               H_\eps(\alpha-\frac{\pi}{2})\right),
\eeq
and $\sigma(\eps):=|\log(\eps)|^{-1}$ has a unique solution $x=[(x_\eps)_\eps]$ 
in $\gs[\R,S^1]$.

Moreover if $\alpha_0\in (-\frac{\pi}{2},\frac{\pi}{2})$, $x_\eps$ has 
the following (continuous) pointwise limit
\beq\label{limit}
 x_\eps(t)\,\to\,x_{\alpha_0}(t)\,:=\,
  \left\{\begin{array}{crcl} 
   e^{-i\frac{\pi}{2}}& -\infty\ < &t&\leq\ -\alpha_0-\frac{\pi}{2}\\
   e^{i(\alpha_0+t)}\qquad& -\alpha_0-\frac{\pi}{2}
                                \ \leq &t&\leq\ -\alpha_0+\frac{\pi}{2}\\
   e^{i\frac{\pi}{2}}& -\alpha_0+\frac{\pi}{2}\ \leq &t&<\ \infty.
  \end{array}\right.
\eeq 
\ethr

\proof 
Existence and uniqueness follows by Theorem \ref{th:mfode} due 
to our assumptions on $\sigma$. 

To prove the statement on the limit we use the following notation:
$x_\eps(t)=e^{i\gamma_\eps(t)}$. 
First note that for all $\alpha_0\in(-\frac{\pi}{2},\frac{\pi}{2})$
we have $\xi_\eps(e^{i\alpha_0})=(e^{i\alpha_0},1)$ for $\eps$ small enough. 
So $e^{i(\alpha_0+t)}$ is a solution as long as $-\frac{\pi}{2}+\sigma(\eps)
<\alpha_0+t=\gamma_\eps(t)<\frac{\pi}{2}-\sigma(\eps)$, that is $-\frac{\pi}{2}+\sigma(\eps)
-\alpha_0<t<\frac{\pi}{2}-\sigma(\eps)-\alpha_0$. Hence $\gamma_\eps(t)\to
\alpha_0+t$ for $-\frac{\pi}{2}-\alpha_0<t<\frac{\pi}{2}-\alpha_0$. 

On the other hand if $t\leq-\frac{\pi}{2}+\sigma(\eps)-\alpha_0$ resp.\
$t\geq\frac{\pi}{2}-\sigma(\eps)-\alpha_0$ then $-\frac{\pi}{2}-\sigma(\eps)
\leq\gamma_\eps(t)\leq-\frac{\pi}{2}+\sigma(\eps)$ resp.\ $\frac{\pi}{2}-\sigma(\eps)
\leq\gamma_\eps(t)\leq\frac{\pi}{2}+\sigma(\eps)$ by the fact that $e^{-i(\frac{\pi}{2}+
\sigma(\eps))}$ resp.\ $e^{i(\frac{\pi}{2}+\sigma(\eps))}$ are equilibrium points and the 
monotonicity of $\gamma_\eps$. Hence the claim follows.
\ep

Note that if $\alpha_0\not\in[-\frac{\pi}{2},\frac{\pi}{2}]$ the
solution equals $e^{i\alpha_0}$ for all times $t$. If $\alpha_0=
\pm\frac{\pi}{2}$ the limit of the solution will in general depend on
the choice of $\rho$. The most ``generic'' choice is (a) to suppose that
$0<\gamma_-\leq H_\eps(0)\leq\gamma_+<1$ for all $\eps$. In this case we obtain
\beq\label{generic}
 x_\eps(t)\to x_{\pm\frac{\pi}{2}}(t)=
 \left\{\begin{array}{crcl}
  e^{-i\frac{\pi}{2}}
   &-\infty\ <&t&\leq\ \mp\frac{\pi}{2} -\frac{\pi}{2}\\
  e^{i(\pm\frac{\pi}{2}+t)}
   &\mp\frac{\pi}{2}-\frac{\pi}{2}\ \leq &t&\leq\ \mp\frac{\pi}{2}+\frac{\pi}{2}\\
  e^{i\frac{\pi}{2}}
   &\mp\frac{\pi}{2}+\frac{\pi}{2}\ \leq &t& <\ \infty\,.
 \end{array}\right.
\eeq
Indeed for $\alpha_0=-\frac{\pi}{2}$ (for $\alpha_0=\frac{\pi}{2}$ 
just adapt the argument accordingly) and $t\leq 0$ we use the same 
arguments as in the last part of the above proof to conclude that 
$-\frac{\pi}{2}-\sigma(\eps)\leq\gamma_\eps(t)\leq-\frac{\pi}{2}$. 
Hence $\gamma_\eps(t)\to-\frac{\pi}{2}$ for $0\leq t$.

To deal with nonnegative $t$ we first observe that
$\mathrm{pr}^2(\xi_\eps(e^{-i\frac{\pi}{2}}))
=H_\eps(0)\geq\gamma_->0$ hence $\dot\gamma_\eps(t)\geq\gamma_-$
for all $t\ge 0$ small enough, i.e., such that
$\gamma_\eps(t)\leq \frac{\pi}{2}-\sigma(\eps)$. So for all such $t$
we obtain $\gamma_\eps(t)\geq\gamma_-\cdot t-\frac{\pi}{2}$. In particular
for $t>\sigma(\eps)/\gamma_-$ we have
$\gamma_\eps(t)\geq-\frac{\pi}{2}+\sigma(\eps)$. So there exists
$t_\eps\leq\sigma(\eps)/\gamma_-$ such that
$\gamma_\eps(t_\eps)=-\frac{\pi}{2}+\sigma(\eps)$. This in turn
implies that for $t_\eps\leq t\leq\pi-2\sigma(\eps)+t_\eps$ the
solution takes the form $\gamma_\eps(t)=
-\frac{\pi}{2}+\sigma(\eps)+(t-t_\eps)$. So
$\gamma_\eps(t)\to-\frac{\pi}{2}+t$ for $0\leq t\leq\pi$.

Finally for $t\geq\pi$ we again use the monotonicity of $\gamma_\eps$ and the fact that 
$e^{i(\frac{\pi}{2}+\sigma(\eps))}$ is an equilibrium point to establish the claim.
\medskip

If we choose (b) $H_\eps(0)=0$ (resp.\ (c) $H_\eps(0)=1$) one sees by adapting the 
above line of arguments
that the limiting solution with initial value $\alpha_0=-\frac{\pi}{2}$ 
($\alpha_0=\frac{\pi}{2}$) will be trapped at $e^{i\alpha_0}$ and equal 
$x_{\alpha_0}$ for $\alpha_0=\frac{\pi}{2}$ ($\alpha_0=-\frac{\pi}{2}$).
However, we still could use cases (b) and (c) in the construction to follow.
In case we drop the assumption $\rho\geq 0$ the limiting behavior can be more complicated
since the solution then may be trapped between different equilibria. 
\medskip

Now we are going to show that the above proposition provides a
counter-example to Marsden's theorem. We consider the flow
$\Phi^\eps(t,e^{i\alpha})=x_\eps(t)$, where $x_\eps$ is the solution
with $x_\eps(0)=e^{i\alpha}$ provided by the Proposition. By Theorem 
\ref{th:mfflow} $\Phi=[(\Phi^\eps)_\eps]$ is in
$\gs[R\times S^1,S^1]$ and has the flow properties (\ref{eq:mfflow2}),
(\ref{eq:mfflow3}). Defining $\Psi=\lim_{\eps\to 0}\Phi^\eps$ conditions
(i) (even with convergence everywhere) and (ii) of Marsden's definition are
satisfied, hence $\zeta$ (given by eq.\ (\ref{zeta})) is a vector
field with measurable flow $\Psi$. So $\Psi$ ought to have the flow
property in the sense of \cite{marsden}, Th.\ 6.2. However,
we have by the second part of Proposition \ref{counterex} and the remark following
its proof (using case (a))
\beq
\Phi^\eps(t,e^{i\alpha})\,\to\,\Psi(t,e^{i\alpha})\,:=\,
\left\{\begin{array}{cl}
    x_{\alpha}(t) &\mbox{if }\ \alpha\in [-\frac{\pi}{2},\frac{\pi}{2}]\\
    e^{i\alpha}&\mbox{if }\ \alpha\not\in
    [-\frac{\pi}{2},\frac{\pi}{2}],
 \end{array}\right.
\eeq 
where $x_\alpha(t)$ denotes the limiting function in (\ref{limit})
with $\alpha=\alpha_0$. Hence 
\beq
 \Psi(-\pi,e^{i\alpha})\,=\,
 \left\{\begin{array}{cl}
  e^{-i\frac{\pi}{2}}&\mbox{if }\ \alpha\in [-\frac{\pi}{2},\frac{\pi}{2}]\\
  e^{i\alpha}&\mbox{otherwise}
 \end{array}\right.
\eeq
and
\beq 
 \Psi(\pi,e^{i\alpha})\,=\,
 \left\{\begin{array}{cl}
  e^{i\frac{\pi}{2}}&\mbox{if }\ \alpha\in [-\frac{\pi}{2},\frac{\pi}{2}]\\ 
  e^{i\alpha}&\mbox{otherwise.}
 \end{array}\right.
\eeq
This in turn implies
\beq
 \Psi(\pi,\Psi(-\pi,e^{i\alpha}))\,=\,
 \left\{\begin{array}{cl}
  e^{i\frac{\pi}{2}}&\mbox{if }\ \alpha\in [-\frac{\pi}{2},\frac{\pi}{2}]\\ 
  e^{i\alpha}&\mbox{otherwise.}
 \end{array}\right.
\eeq
So $\Psi_\pi\circ\Psi_{-\pi}\not={\mathrm id}$ for all 
$e^{i\alpha}$ with $\alpha\in[-\frac{\pi}{2},\frac{\pi}{2}]$ 
contradicting the assertion of Theorem 6.2 in \cite{marsden}. 
Note that if we choose cases (b) or (c) above we obtain similar results; only
the range of $\alpha$ changes form the closed interval to the half open
interval $(-\frac{\pi}{2},\frac{\pi}{2}]$ resp.\ $[-\frac{\pi}{2},\frac{\pi}{2})$
and again the flow property fails to hold for $e^{i\alpha}$ in 
a set of positive measure.
\medskip

Since the approach in \cite{marsden} is built upon pointwise
convergence almost everywhere of the regularizing flows 
its failure motivates the study of different notions of convergence for 
generalized functions taking values in a manifold to allow for a corrected 
version of the Theorem. We do so in the following section.

\section{Notions of Association}\label{sec:ass}

In all variants of spaces of Colombeau generalized functions taking values in a
linear space compatibility with respect to the distributional setting
is affected through the notion of association. We call $u\in\gs(X)$
associated with zero, $u\approx 0$, if one (hence every) representative 
$u_\eps$ converges to zero weakly (cf.\ also Definition \ref{def:association} 
(v) below). The assignment $u\approx v 
:\Leftrightarrow u_\eps-v_\eps\approx 0$ gives rise to an equivalence
relation on $\gs(X)$ and a linear quotient space $\gs(X)/\approx$, 
generalizing distributional equality to the level of $\gs(X)$. Moreover
if $\lim_{\eps\to 0}\int_Xu_\eps\nu=\langle\omega,\nu\rangle$ for some 
distribution $\omega$ and every compactly supported one-density $\nu$ 
we write $u\approx \omega$ and call $\omega$ the distributional shadow 
of $u\in\gs(X)$.
 
In this section we are going to introduce a number of notions of association in the
space $\gs[X,Y]$ (cf. also \cite{gfvm2}, Sec.\ 6) and clarify their respective 
interrelations. 

\bd\label{def:association}
 Let $u=[(u_\eps)_\eps],v=[(v_\eps)_\eps]\in\gs[X,Y]$, and let $h$
 be a Riemannian metric on $X$
 with distance function $d_h$. 
 \begin{itemize}
  \item[(i)] $u$ is called {\em zero-associated} (cf.\ \cite{gfvm2}, Def 6.1) with $v$,
   \beq\label{ass0}
    u\approx_0 v\ :\Leftrightarrow \sup_{p\in K}d_h(u_\eps(p),v_\eps(p))\to 0\quad
    \forall K\comp X.
   \eeq
  \item[(ii)] $u$ is called {\em pointwise-associated (pw-associated)} with $v$,
   \beq\label{asspw}
    u\approx_{\mathrm pw} v\ :\Leftrightarrow d_h(u_\eps(p),v_\eps(p))\to 0\quad 
    \forall p\in X.
   \eeq
 \item[(iii)] $u$ is called {\em pointwise-associated almost everywhere 
 (pwae-associated)} 
   with $v$,
   \beq\label{asspwae}
    u\approx_{\mathrm pwae} v\ :\Leftrightarrow d_h(u_\eps(p),v_\eps(p))\to 0\ 
    \mbox{ for almost all }\ p\in X.
   \eeq
 \item[(iv)]  $u$ is called {\em model-associated} with $v$, 
   \beq\label{assmodel}
    u\approx_{\small{\mathcal M}} 
    v\ :\Leftrightarrow f\circ u_\eps-f\circ v_\eps\to 0\ \mbox{ in }\
    \D'(X)\quad \forall f\in\CC^\infty(Y).
   \eeq
 That is, $f\circ u\approx f\circ v$ in $\gs(X)$ for all $f\in
 \CC^\infty(Y)$ (see (v) below). 
 \item[(v)] If $Y=\R^n$ then $u$ is called {\em associated} with $v$,
   \[ u\approx v:\Leftrightarrow u_\eps-v_\eps\to 0\ \mbox{ in }\ \D'(X)^n.
   \]
 \end{itemize}
\ed
It is straightforward to check that notions (i)--(iii) are independent 
of the Riemannian metric $h$ employed (cf.\  \cite{book}, Lemma 3.2.4) 
and that all of the above definitions are independent of the representatives 
chosen for $u$ and $v$.

\bthm \label{th:assoc}
 We have the following chain of implications: 
 \begin{center}
   (i) $\Rightarrow$ (ii) $\Rightarrow$ (iii)  $\Rightarrow$ (iv) $\Rightarrow$ (v),
 \end{center}
 where the last implication holds in case $Y=\R^n$.  None of the above implications can be reversed.
\ethm

\proof
 (i) $\Rightarrow$ (ii) $\Rightarrow$ (iii) is clear as well as (iv) $\Rightarrow$ (v).\\
 (iii) $\Rightarrow$ (iv): Let $f\in\CC^\infty(Y)$ and let $(u_\eps)_\eps$ resp.\ $(v_\eps)_\eps$
  be  representatives of $u$ resp.\ $v$. Then $|f\circ u_\eps(p)-f\circ v_\eps(p)|\to 0$ a.e.\ 
  and is bounded uniformly in $\eps$ on compact sets by the c-boundedness of 
  $u_\eps$ and $v_\eps$.
  Hence by dominated convergence $\int\left(f\circ u_\eps(p)-f\circ v_\eps(p)\right)\mu(p)
  \to 0$ for all compactly supported one-densities $\mu$ on $X$.\\
 (v) $\not\Rightarrow$ (iv): Set $u_\eps=\sin\left(\frac{x}{\eps}\right)$ and $v_\eps=0$.
  Then $u\approx v$ in $\gs(\R)$ but for $f(x)=x^2$ we have $f\circ u_\eps\to 1/2$ pointwise
  hence $f\circ u_\eps-f\circ v_\eps\not\to 0$ in $\D'(\R)$.\\
 (iv) $\not\Rightarrow$ (iii): Let $\rho_0\in\D(\R)$,
  $\rho_0(\R)\subseteq[0,1]$, $\supp(\rho_0)
  \subseteq[-1,1]$, $\int\rho_0=1$ and $\rho_0(0)=1$. Furthermore, set
  \[ \rho(x):=\sum\limits_{n=-\infty}^{n=\infty}\rho_0(2^{|n|}(x-n))
  \]
  and $\rho_\eps(x):=\rho(x/\eps)$. Then there is no $x\in\R$ such that 
  $\rho_\eps(x)\to 0$. Indeed,
 $\forall x\ \forall n\ \exists \eps_n$, $\eps_n\to 0$ such that
  $\rho_{\eps_n}(x)=1$. On the other hand $\rho_\eps\approx_{\mathcal M} 0$ since 
  for any $f\in\CC^\infty(\R)$ and any test function $\varphi$ we have
  \beas
   |\int(f\circ \rho_\eps(x)) \varphi(x)\, dx|
   &\leq& \|\nabla f\|_{\infty,[0,1]}\int|\rho_\eps(x)|\,|\varphi(x)|\,dx\\
   &=& \|\nabla f\|_{\infty,[0,1]}\int\rho\left(\frac{x}{\eps}\right)|\varphi(x)|\, dx\\
   &\leq& \eps \|\nabla f\|_{\infty,[0,1]}\, \|\varphi\|_\infty\, \|\rho\|_1\\
   &\to& 0\,.
  \eeas
 (iii) $\not\Rightarrow$ (ii) $\not\Rightarrow$ (i) is clear.
\ep

\bd 
 Let $u=[(u_\eps)_\eps]$ in $\gs[X,Y]$ and $v:X\to Y$ a map. 
 $v$ is called a {\em shadow} of $u$ in the sense of zero-, pw-,
 pwae-, resp.\ model-association (or, for short, zero-, pw-, pwae-, resp.\ 
 model-associated with $u$) if (\ref{ass0}), (\ref{asspw}), (\ref{asspwae}), 
 resp.\ (\ref{assmodel}) holds with $v$ replacing $v_\eps$.
\ed

Of course the hierarchy of Theorem \ref{th:assoc} carries over to shadows of the
types introduced above. 
In the following we shall also need a notion which encodes information
on the order of convergence with respect to $\eps$.

\bd
 Let $u\in\gs[X,Y]$. 
 \begin{itemize}
 \item[(i)] $u$ is called {\em fast-associated} with $v\in\gs[X,Y]$, 
 \[
  u\approx_{\mathrm f}v\ :\Leftrightarrow d_h(u_\eps(p),v_\eps(p))=O(\eps^m)
  \ \forall p\in X \ \forall \eps \in \N
 \]
 for one (hence every) representative $(u_\eps)_\eps$ of $u$ and one 
 (hence every) representative $(v_\eps)_\eps$ of $v$ where again
 $d_h$ denotes the Riemannian distance with respect to any  
 Riemannian metric $h$ on $Y$.
\item[(ii)] $u$ is called {\em fast-associated} with $v:X\to Y$, 
 \[
  u\approx_{\mathrm f}v\ :\Leftrightarrow d_h(u_\eps(p),v(p))=O(\eps^m)
  \ \forall p\in X \ \forall \eps \in \N
 \]
 for one (hence every) representative $(u_\eps)_\eps$ of $u$.
 \end{itemize}

\ed

Note that for $u\,, v\in\gs[X,Y]$, $u\approx_{\mathrm f}v$ if and only if
the generalized number $d_h(u(p),v(p))$ is $0$ for each $p\in X$ and one
(each) Riemannian metric $h$, i.e., if and only if $u(p) = v(p)$ for all $p\in X$.
This notion is {\em strictly} weaker than equality in $\gs[X,Y]$ (cf.\
\cite{point}). $\approx_{\mathrm f}$ implies $\approx_{{\mathrm pw}}$ and the
converse implication is clearly wrong. However, there is no relation between
$\approx_{\mathrm f}$ and $\approx_0$.

\section{Limiting flows}\label{sec:lim}

Having introduced a number of notions of association in the previous
section we now have the tools at hand to analyze the following
question:
Let $\xi\in\gs^1_0(X)$ be $\gs$-complete
with $\Phi=[(\Phi^\eps)_\eps]$ its (unique) flow in 
$\gs[\R\times X,X]$. If $\Phi$ admits a shadow $\Psi$ in the sense of one of the 
notions introduced above, does this imply that $\Psi$ has the 
flow property? We are going to answer this question in the following but first 
turn to some preliminaries.
\medskip

\bd 
 We say a function $u\in\gs[X,Y]$ is of  {\em locally bounded derivative} if
 for all $K\comp X$ and for one (hence every) pair of Riemannian metrics $g$ on 
 $X$ resp.\ $h$ on $Y$ there exists $C>0$ and $\eps_0>0$ such that
 \[
  \sup\limits_{p\in K}\|T_pu_\eps\|_{g,h}\leq C\qquad \forall \eps\leq\eps_0
 \]
 for one (hence every) representative $(u_\eps)_\eps$ of $u$. Here $\|T_pf\|_{g,h}$
denotes the norm of the linear map $T_pf:\ (T_pX,\|\quad\|_g)\to(T_{f(p)}Y,\|\quad\|_h)$
(cf.\ \cite{book}, 3.2.54).
\ed

\blem\label{th:lemflow}
 Let $(X,h)$ be a complete Riemann manifold, let $\xi\in\gs^1_0(X)$ 
 satisfy the assumptions of Theorem \ref{th:mfode} and denote by
 $\Phi = [(\Phi^\eps)_\eps]\in \gs[\R\times X,X]$ the generalized flow of $\xi$.
   If $P\xi$ is locally bounded for all differential operators 
   $P\in{\cal P}(X,TX)$ of first order then $\Phi(t,.)$ is of locally 
   bounded derivative. 
\ethr

\proof
 For each $\eps$ let $\Phi^\eps$ be the complete flow corresponding to some 
 globally bounded (w.r.t. $h$) representative $\xi_\eps$ of $\xi$. 
 Then each $T\Phi^\eps(t,.)$ satisfies the following ODE 
 \beas
  \frac{d}{dt}T\Phi^\eps(t,\,.\,)&=&T\xi(\Phi^\eps(t,\,.\,))\ T\Phi^\eps(t,\,.\,)\\ 
  T\Phi^\eps(0,\,.\,)&=&\mathrm{id}.
 \eeas
 As we only need to consider $(t,p)$ varying in a compact subset of
 $\R\times X$ we may employ the same argument as in the
 uniqueness part of the proof of Theorem \ref{th:mfflow} to
 successively estimate over pieces of the integral curves
 $s\mapsto\Phi^\eps(s,p)$ of $\xi_\eps$,  each contained in a 
 single chart. Therefore we may work locally to obtain
 \[ 
  D_x\Phi^\eps(t,x)=x+\int\limits_0^t D_x\xi_\eps(\Phi^\eps(s,x))\ D_x\Phi^\eps(s,x)\, ds
 \]
 and by Gronwall's inequality 
 \[ 
  \|D\Phi^\eps(t,x)\|\leq C\ \exp\left(\int\limits_0^t \|D_x\xi_\eps(\Phi^\eps(s,x))\|\,ds\right).
 \]
 \hspace*{0pt} From this and the c-boundedness of the flow the assertion 
 follows.
\ep

\bthm\label{th:limflow}
 \begin{itemize}
  \item [(i)] Let $\xi\in\gs^1_0(X)$ satisfy the assumptions of Lemma \ref{th:lemflow}.
    If for all $t\in \R$, $\Phi(t,\,.\,)\approx_{\mathrm f}\Psi(t,\,.\,)$ 
    then $\Psi$ has the flow property.
  \item [(ii)] Let $\xi\in\gs^1_0(X)$ be a $\gs$-complete generalized vector field with flow
    $\Phi$. If for each $t\in\R$, $\Phi(t,\,.\,)\approx_0\Psi(t,\,.\,)$ then $\Psi$ has the flow property.
 \end{itemize}
\ethm

\proof 
 (i) On the complete Riemannian manifold ($X,h)$ we have to show that
 \[
  d_h(\Psi(s+t,p),\Psi(s,\Psi(t,p)))=0\,.
 \]
 for all $s,t\in\R$ and for all $p\in X$. Since for all $s,t,p$
 \[
  \Psi(s+t,p)=\lim_{\eps\to 0}\Phi^\eps(s+t,p)=\lim_{\eps\to 0}\Phi^\eps(s,\Phi^\eps(t,p))
 \]
 it suffices to show that
 \[
  d_h(\Phi^\eps(s,\Phi^\eps(t,p)),\Psi(s,\Psi(t,p)))\to 0
 \]
 as $\eps\to 0$. We introduce the following splitting 
 \beq\label{eq:split}
  d_h(\Phi^\eps_s(\Phi^\eps_t(p)),\Psi_s(\Psi_t(p)))
  \,\leq\,d_h(\Phi^\eps_s(\Phi^\eps_t(p)),\Phi^\eps_s(\Psi_t(p)))
   +d_h(\Phi^\eps_s(\Psi_t(p)),\Psi_s(\Psi_t(p))).
 \eeq
 Here the second term converges to zero since $\Phi^\eps(s,q)\to\Psi(s,q)$ pointwise
 and we are left with the first term. On a complete Riemannian manifold
 any two points can be joined by a minimizing geodesic segment. So we choose
 such segments $\ga_\eps:[0,b_\eps]\to X$ with $\ga_\eps(0)=\Psi(t,p)$ and
 $\ga_\eps(b_\eps)=\Phi^\eps(t,p)$. Then by the $c$-boundedness of
 $\Phi^\eps$ and the Hopf-Rinow theorem we may choose $K\comp X$ and
 $\eps_0>0$ such  that each $\gamma_\eps$ stays entirely within $K$ 
 for $\eps<\eps_0$. Thus
 \bea\label{xxx}
   d_h(\Phi^\eps(s,\Phi^\eps(t,p)),\Phi^\eps(s,\Psi(t,p)))
    &\leq&\int\limits_0^{b_\eps}\|(\Phi^\eps(s,.)\circ\ga_\eps)'(\lambda)\|_h\,d\lambda\nn\\
   &\leq&\int\limits_0^{b_\eps}\|T_{\ga_\eps(\lambda)}\Phi^\eps(s,.)\|_{h,h}\,
    \|\ga_\eps'(\lambda)\|_h\,d\lambda\nn\\
   &\leq&\sup\limits_{q\in K}\|T_q\Phi^\eps(s,.)\|_{h,h}\,
    \int\limits_0^{b_\eps} \|\ga_\eps'(\lambda)\|_h\,d\lambda\\
   &\leq& \frac{C}{\eps^N}\,\,d_h(\Phi^\eps(s,p)),\Psi(s,p)),\nn
 \eea
 for some $C,\, N > 0$. Since $\Phi(t,p)  \approx_{\mathrm f}\Psi(t,p)$, 
 $d_h(\Phi^\eps(s,\Phi^\eps(t,p)), \Phi^\eps(t,\Psi(s,p)))$ converges
 to zero for $\eps\to 0$, as desired.

 (ii) In this case we use the splitting 
 \beq\label{eq:split2}
  d_h(\Phi^\eps_s(\Phi^\eps_t(p)),\Psi_s(\Psi_t(p)))
  \,\leq\,d_h(\Phi^\eps_s(\Phi^\eps_t(p)),\Psi_s(\Phi^\eps_t(p)))
   +d_h(\Psi_s(\Phi^\eps_t(p)),\Psi_s(\Psi_t(p)))
 \eeq
 Here the first term converges to $0$ since $\Phi^\eps(s,\,.\,) \to \Psi(s,\,.\,)$
 locally uniformly and $(\Phi^\eps)_\eps$ is c-bounded and the second
 one since $\Psi_s$ is necessarily continuous.
\ep

As was already mentioned in the proof, in case (ii) the limiting flow
$\Psi$ necessarily is continuous in $p$. On the other hand, we will
give an explicit example of a discontinuous limiting flow below.
First, however, we turn to another set of assumptions
guaranteeing a continuous limiting flow.

\bc
 Let  $\xi\in\gs^1_0(X)$ satisfy the assumptions of \ref{th:lemflow} (ii).
 For each $t\in \R$, let $\Phi(t,\,.\,)\approx_{{\mathrm
 pw}}\Psi(t,\,.\,)$. Then $\Psi$ is continuous in $p$ and  has the flow property.
\ethr

\proof 
 Let $\Phi = [(\Phi^\eps)_\eps]$. By Lemma \ref{th:lemflow} (ii) it
 follows that $(\Phi^\eps)_\eps$ is locally uniformly equicontinuous.
 But then by the theorem of Arzela-Ascoli $\Phi^\eps$ in fact
 converges locally uniformly to $\Psi$, i.e.,
 $\Phi(t,\,.\,)\approx_0\Psi(t,\,.\,)$ for all $t$. Hence the claim
 follows from Theorem \ref{th:limflow} (ii).
\ep

\bex
 Let $X=T^2 = S^1\times S^1$ and $\xi=[(\xi)_\eps] = \gs^1_0(X)$ be given by
 \[
  \xi_\eps(e^{i\alpha},e^{i\beta})=(e^{i\alpha},e^{i\beta};1,1-\rho_{\sigma(\eps)}(\alpha)).
 \]
 Although  $\xi$ does not satisfy the boundedness assumption of Theorem
 \ref{th:mfflow}, we may nevertheless establish its $\gs$-completeness
 as follows. First, since $X$ is compact, each $\xi_\eps$ possesses a
 global flow $\Phi^\eps$. $(\Phi^\eps)_\eps$ is moderate by the proof
 of Theorem \ref{th:mfflow} and that $\Phi:=[(\Phi^\eps)_\eps]$ is
 indeed the unique flow of $\xi$ follows readily by choosing
 appropriate charts on $X$ and applying the local case (Theorem
 \ref{th:locflow}).   

 Moreover, $\Phi$ has a discontinuous pointwise limit $\Psi$, namely
 \beq
  \Phi^\eps(t;e^{i\al},e^{i\beta})
  =\left(\begin{array}{l}
   e^{i(\alpha +t)}\\e^{i(\beta+t-\int\limits_\alpha^{\alpha+t}
   \rho_{\sigma(\eps)}(\gamma)\,d\gamma)}
  \end{array}\right)
  \,\to\,
  \left(\begin{array}{l}
   e^{i(\al+t)}\\e^{i(\beta+t-H(\alpha+t)+H(\alpha))}
  \end{array}\right),
 \eeq
  (in fact we even have $\Phi^\eps(t,\,.\,)\approx_{\mathrm f} \Psi(t,\,.\,)$ for all $t$) 
  which by  direct verification  satisfies the flow property  
  $\Psi_{s+t}=\Psi_s\circ\Psi_t$ for all $s.\, t\in\R$. 
\ethr

To conclude this section let us summarize the above results in the following way. From the 
counter-example in Section \ref{distsect} we know that $\Phi\approx_{\mathrm{pw}}\Psi$ 
does not imply the flow property of $\Psi$. Neither does convergence
of $\Phi^\eps(t,.)$ to $\Psi(t,.)$ locally in $L^p$ for any $1\leq
p<\infty$ secure the flow property of $\Psi$, as can also be seen from
the example given in Section \ref{distsect}. On the other hand, convergence 
locally in $L^\infty$ (i.e., (ii) of Theorem \ref{th:limflow}) implies the flow 
property of $\Psi$ while at the same time forcing the limiting flow to
be continuous. Hence 
the above example lies precisely in the gap which allows for a discontinuous limiting {\em flow}.

\end{document}